\documentstyle{amsart}

\title[Algebraic Systems...] {Algebraic Systems with\\ Lipschitz Perturbations}

\author{Giovanni Molica Bisci and Du\v{s}an Repov\v{s}}
\address[G. Molica Bisci]{Dipartimento P.A.U., Universit\`a  degli
Studi Mediterranea di Reggio Calabria, Salita Melissari - Feo di
Vito, 89124 Reggio Calabria, Italy} \email{gmolica@@unirc.it}

\address[D. Repov\v{s}]{Faculty of Education, and Faculty of Mathematics and Physics\\ University of Ljubljana, POB 2964, Ljubljana, Slovenia 1001}
\email{dusan.repovs@@guest.arnes.si}



{}
\newtheorem{theorem}{Theorem}[section]
\newtheorem{lemma}{Lemma}[section]
{}

\newtheorem{remark}{Remark}[section]

\def\trace{\mathop{\rm Tr}(A)}

\newcommand{\erre}{\mbox{\normalshape I\!R}}
\newcommand{\enne}{\mbox{\normalshape I\!N}}
\begin{document}

\begin{abstract}
 By using variational methods, the existence of infinitely many solutions for a nonlinear algebraic system with a parameter is established in presence of a perturbed Lipschitz term. Our goal was achieved requiring an appropriate behavior of the nonlinear term $f$, either at zero or at infinity, without symmetry conditions.
\end{abstract}
\maketitle
\section{Introduction}
In many cases a problem in a continuous framework can be handled
by using a suitable method from discrete mathematics, and
conversely. The modeling/simulation of certain
nonlinear problems from economics, biological neural networks,
optimal control and others, enforced in a natural manner a rapid
development of the theory of discrete equations.\par
\indent In this paper, motivated by this large interest, we study the following algebraic system
$$
Au=\lambda f(u)+h(u),\eqno\displaystyle{(S_{A,\lambda}^{f,h})}
$$
 in which $u=(u_1,...,u_n)^{t}\in\erre^n$ is a column vector, $A=(a_{ij})_{n\times n}$ is a positive-definite matrix, $f(u):=(f_1(u_1),...,f_n(u_n))^t$, where the functions $f_k:\erre\rightarrow\erre$ are assumed to be continuous
for every $k\in {\mathbb{Z}}[1,n]:=\{1,2,...,n\}$,
and $\lambda$ is a positive parameter.\par
\indent
 Moreover, $$h(u):=(h_1(u_1),...,h_n(u_n))^t,$$ where, for every $k\in {\mathbb{Z}}[1,n]$, the functions $h_k:\mathbb{R} \to \mathbb{R} $ are
Lipschitz
continuous with constants $L_k\geq0$, that is:
$$|h_k(t_1)-h_k(t_2)|\leq L_k|t_1-t_2|,$$
for every $t_1,t_2\in \mathbb{R}$, and $h_k(0)=0$.
\par
A large number of discrete problems can be formulated as special cases of the non-perturbed ($h=0$) algebraic system, namely $(S_{A,\lambda}^{f})$; see, for instance, the papers \cite{YZ, Z,ZB, ZC, ZF} and references therein.\par
 We also point out that the special case
$$
A:=\left(\begin{array}{ccccc}
  2 & -1 & 0 & ... & 0 \\
  -1 & 2 & -1 & ... & 0 \\
     &  & \ddots &  &  \\
  0 & ... & -1 & 2 & -1 \\
  0 & ... & 0 & -1 & 2
\end{array}\right)_{n\times n},
$$
 has been considered in
order to study the existence of nontrivial solutions of nonlinear
second-order difference equations \cite{KMRT,KMR,MRT}. Further, general references on difference equations and their applications can be found e.g. in \cite{A,kelly}.\par\par
 Here, by using variational methods, under the key assumption that
  $$
  \displaystyle L:=\max_{k\in {\mathbb{Z}}[1,n]}L_k<\lambda_1,
  $$
   where $\lambda_1$ is the first eigenvalue of the matrix $A$, we determine open intervals of positive parameters such that problem $(S_{A,\lambda}^{f,h})$ admits either an unbounded sequence of solutions, provided that the nonlinearity $f$ has a suitable behaviour at infinity (Theorem \ref{Ma}), or a sequence of pairwise distinct solutions that converges to zero, if a similar behaviour occurs at zero (see Theorem \ref{Main2}).\par
 Our main tool is a recent critical point result obtained by Ricceri and recalled here in a convenient form (see Theorem \ref{two}).\par

 A special case of our results reads as follows (see Remark \ref{din}).
\begin{theorem}\label{intro}Let $z: \erre \rightarrow \erre$ be a nonnegative and continuous function. Assume that

\[
\displaystyle{\liminf_{t\rightarrow +\infty}}\frac{\displaystyle\int_{0}^{t}z(\xi)d \xi}{t^2}=0,\,\,\,\,\,\,\,\,\,\limsup_{t\rightarrow+\infty}\frac{\displaystyle\int_{0}^{t}z(\xi)d \xi}{t^2}=+\infty.
\]

\noindent Then, for each $\lambda>0,$ and for every Lipschitz function $h:\erre\rightarrow \erre$ with sufficiently small constant $L_h,$
the following discrete problem
$$
[u(i+1,j)-2u(i,j)+u(i-1,j)]+[u(i,j+1)-2u(i,j)+u(i,j-1)]
$$
$$
+\lambda z(u(i,j))+h(u(i,j))=0,\,\,\,\,\,\,\,\,\forall\,(i,j)\in {\mathbb{Z}}[1,m]\times {\mathbb{Z}}[1,n]
$$
with boundary conditions
$$
u(i,0)=u(i,n+1)=0,\,\,\,\,\,\,\forall\, i\in {\mathbb{Z}}[1,m],
$$
$$
u(0,j)=u(m+1,j)=0,\,\,\,\,\,\,\forall\, j\in {\mathbb{Z}}[1,n],
$$
admits an unbounded sequence of solutions.
\end{theorem}

 Finally, for completeness, we just mention here that there is a vast literature on nonlinear difference equations based on fixed point and upper and lower solution methods (see \cite{BoMawhin, Henderson}). For related topics see the works \cite{A3,A4,A2,A1}. For a complete and exhaustive overview on variational methods we refer the reader to the
monographs \cite{KRV,MR}.

\section{Abstract Setting}

Let $(X,\|\cdot\|)$ be a finite-dimensional Banach space and let
$J_\lambda:X\rightarrow \erre$ be a function satisfying the
following structure hypothesis:
\begin{itemize}
\item[$(\Lambda)$] for all $u\in X$, \textit{$J_\lambda(u):=\Phi(u)-\lambda\Psi(u)$
 where $\Phi, \Psi:X\rightarrow \erre$ are
two functions of class $C^1$ on $X$ with $\Phi$ coercive, i.e.
$\lim_{\|u\|\rightarrow \infty}\Phi(u)=+\infty$, and $\lambda$ is
a real positive parameter.}
\end{itemize}
\noindent Moreover, provided that $r>\inf_{X}\Phi$, put
\[
\varphi(r):=\inf_{u \in \Phi^{-1}(\left]-\infty, r\right[)}
\frac{\left(\displaystyle\sup_{v \in {\Phi^{-1}\left(\left]-\infty,
r \right[ \right)}}\Psi(v)\right)-\Psi(u)}{r - \Phi(u)},
\]
\noindent and
$$
\gamma:=\liminf_{r\rightarrow+\infty}\varphi(r),\,\,\,\,\,\,\
\delta:=\liminf_{r\rightarrow (\inf_X\Phi)^{+}}\varphi(r).
$$
\noindent Clearly, $\gamma \ge 0$ and $\delta \ge 0$. When $\gamma =0$ (or $\delta = 0$), in the sequel, we agree to read $1/\gamma$ (or
$1/\delta$) as $+\infty$.\par
\begin{theorem} \label{two}
Assuming that the condition $(\Lambda)$ holds, one has
\begin{itemize}
\item [$(\textrm{a})$] If $\gamma<+\infty$ then, for each
$\lambda\in \left]0,{1}/{\gamma}\right[$, the following alternative
holds$:$\\
either
\begin{itemize}
 \item[$(\textrm{a}_1)$] $J_{\lambda}$ possesses a global minimum,
\end{itemize}
or
\begin{itemize}
\item[$(\textrm{a}_2)$] there is a sequence $\{u_m\}$ of critical points
  $($local minima$)$ of $J_\lambda$ such that
  $\lim_{m\rightarrow \infty}\Phi(u_m)=+\infty$.
\end{itemize}
\end{itemize}
\begin{itemize}
\item [$(\textrm{b})$] If $\delta<+\infty$ then, for each
$\lambda\in \left]0,{1}/{\delta}\right[$, the following alternative
holds$:$\\
either
\begin{itemize}
 \item[$(\textrm{b}_1)$] there is a global minimum of $\Phi$ which is a local minimum of $J_\lambda$,
\end{itemize}
or
\begin{itemize}
\item[$(\textrm{b}_2)$]there is a sequence $\{u_m\}$ of pairwise distinct critical points
  $($local minima$)$ of $J_\lambda$, with $\lim_{m \to \infty}\Phi(u_m)=\inf_X\Phi$, which
  converges to a global minimum of $\Phi$.
\end{itemize}
\end{itemize}
\end{theorem}
\begin{remark}\rm{Theorem \ref{two}
is the 
finite-dimensional version of the quoted multiplicity result of Ricceri from \cite{Ricceri}.}
\end{remark}

\indent As ambient space $X$, consider the $n$-dimensional Banach space $\erre^n$
endowed by the norm
$$
\|u\|:=\Big(\sum_{k=1}^{n}u_k^2\Big)^{1/2}.
$$
\indent Set ${\mathfrak{X}}_{n}$ be the class of all symmetric and positive-definite matrices of order $n$. Further, we denote by $\lambda_1,...,\lambda_n$ the eigenvalues of $A,$ ordered as follows $0<\lambda_1\leq...\leq \lambda_n$.\par
It is well-known that if $A\in {\mathfrak{X}}_{n}$, then for every $u\in X$, one has

\begin{equation}\label{immersione2}
\lambda_1\|u\|^2\leq u^{t}Au\leq \lambda_n\|u\|^2,
\end{equation}

\noindent and

\begin{equation}\label{immersione}
\|u\|_\infty\leq \frac{1}{\sqrt{\lambda_1}}(u^{t}Au)^{1/2},
\end{equation}
\noindent where $ \|u\|_{\infty}:=\displaystyle\max_{k\in
{\mathbb{Z}}[1,n]}|u_k| $.

\indent Set
\begin{equation} \label{newfunction}
\Phi(u):=\frac{u^{t}Au}{2}-\sum_{k=1}^{n}H_k(u_k),
\end{equation}
and
\begin{equation} \label{functions}
\Psi(u):=\sum_{k=1}^{n}F_k(u_k), \quad \quad J_\lambda(u):=\Phi(u)-\lambda\Psi(u),
\end{equation}
\noindent for every $u\in X$, where $H_k(t):=\displaystyle\int_{0}^{t}h_k(\xi)d \xi$ and
$F_k(t):=\displaystyle\int_{0}^{t}f_k(\xi)d \xi$, for every $(k,t)
\in {\mathbb{Z}}[1,n]\times \erre$.\par \indent Standard arguments show that
$J_\lambda \in C^1(X,\erre),$ as well as that critical points of
$J_\lambda$ are exactly the solutions of problem
$(S_{A,\lambda}^{f,h})$; see, for instance, the paper \cite{YZ2}.\par

\begin{lemma}\label{OggiH}
Set
\begin{equation} \label{LLL}
\displaystyle L:=\max_{k\in {\mathbb{Z}}[1,n]}L_k<\lambda_1.
\end{equation}
Then the functional $\Phi$ is coercive.
\end{lemma}

\begin{pf}
Bearing in mind \eqref{immersione2}, since $h_k$ is a
Lipschitz continuous function (for every $k\in {\mathbb{Z}}[1,n]$) with constant
$L_k\geq0$ and $h_k(0)=0$, we have
\begin{eqnarray}
\Phi(u)&\geq&\frac{\lambda_1}{2}\|u\|^2-\sum_{k=1}^n|H_k(u_k)|\geq \frac{1}{2}\|u\|^2-\sum_{k=1}^n\left(\int_0^{u_k}|h_k(t)|dt\right)\nonumber\\
 &\geq& \frac{\lambda_1}{2}\|u\|^2-L\sum_{k=1}^n\int_0^{u_k}|t|dt=\frac{1}{2}\|u\|^2-\frac{L}{2}\sum_{k=1}^n u_k^2\nonumber\\
&=&\left(\frac{\lambda_1-L}{2}\right)\|u\|^2\nonumber.
\end{eqnarray}
Hence, by \eqref{LLL}, the above relation implies that the functional $\Phi$ is coercive.
\end{pf}

\section{Main results}

Set
\[
\displaystyle{A_\infty:=\liminf_{t\rightarrow +\infty}}\frac{\displaystyle{\sum_{k=1}^{n}}\displaystyle{\max_{|\xi|\leq
t}F_k(\xi)}}{t^2}, \quad  \quad\textrm{ and} \quad  \quad B^{\infty}:=\limsup_{t\rightarrow+\infty}\displaystyle{\frac{\displaystyle\sum_{k=1}^{n}F_k(t)}{t^2}}.
\]

From now on we shall assume that the functions $h_k:\mathbb{R} \to \mathbb{R}$, for every $k\in {\mathbb{Z}}[1,n]$, are
Lipschitz
continuous with constants $L_k\geq0$ such that condition \eqref{LLL} holds.
\begin{theorem}\label{Ma}
Let $A\in {\mathfrak{X}}_n$ and assume that the following inequality
holds
\begin{itemize}
    \item [$(\textrm{h}_\infty^L)$] \qquad \qquad $A_\infty <\displaystyle \frac{\lambda_1-L}{\trace +2\sum_{i<j}a_{ij}+nL}B^\infty.$
\end{itemize}

\noindent Then, for each $$\lambda \in \displaystyle\left]
\frac{\trace +2\sum_{i<j}a_{ij}+nL}{2B^\infty},
\frac{\lambda_1-L}{2{{A}}_\infty}\right[,$$ problem
$(S_{A,\lambda}^{f,h})$ admits an unbounded sequence of solutions.
\end{theorem}

\begin{pf}
Fix $\lambda$ as in the assertion of the theorem and put
$\Phi$, $\Psi$, $J_\lambda$ as
in \eqref{newfunction} and \eqref{functions}.
Since the critical points of $J_\lambda$ are the solutions of problem $(S_{A,\lambda}^{f,h})$, our aim is to apply Theorem \ref{two} part $(\rm{a})$ to function $J_\lambda$.
Clearly $(\Lambda)$ holds.\par
 Therefore,  our conclusion 
follows provided that $\gamma<+\infty$  as well as that $J_\lambda$ turns out to be unbounded from below. To this end, let $\{c_m\}$ be a real sequence such
that $\displaystyle \lim_{m\rightarrow \infty}c_m=+\infty$ and
\[
\displaystyle{\lim_{m\rightarrow \infty}}\frac{\displaystyle{\sum_{k=1}^{n}}\displaystyle{\max_{|\xi|\leq
c_m}F_k(\xi)}}{c_m^2}=A_\infty,
\]

\indent Write \[r_m:=\frac{\lambda_1-L}{2}c_m^2,\]
for every $m\in\enne$.\par
 Since,
owing to $(\ref{immersione})$, it follows that
$$
\{v\in X:v^tAv<2r_m\}\subset \{v\in X:|v_k|\leq c_m\,\,\forall\, k\in {\mathbb{Z}}[1,n]\},
$$
and we obtain
$$
\varphi(r_m)\leq \frac{\displaystyle\sup_{v^{t}Av<
2r_m}\displaystyle\sum_{k=1}^{n}F_k(v_k)}{r_m}
\leq
\frac{\displaystyle\sum_{k=1}^{n}\max_{|t|\leq c_m}F_k(t)}{r_m}=\frac{2}{\lambda_1-L}\frac{\displaystyle\sum_{k=1}^{n}\max_{|t|\leq c_m}F_k(t)}{c_m^2}.
$$
\indent Hence, it follows that
$$
\gamma\leq \lim_{m\rightarrow \infty}\varphi(r_m)\leq\displaystyle \frac{2}{\lambda_1-L}A_\infty < \frac{1}{\lambda}<+\infty.
$$
\indent Now, we verify that $J_\lambda$ is unbounded from below.
First, assume that $B^\infty=+\infty$. Accordingly, fix such $M$ that $$M>\displaystyle\frac{\trace +2\sum_{i<j}a_{ij}+nL}{2\lambda}$$ and let $\{b_m\}$ be
a sequence of positive numbers,
with $\displaystyle\lim_{m\rightarrow \infty}b_m=+\infty$, such that
\[
\sum_{k=1}^{n}F_k(b_m)>Mb^2_m, \quad \quad  (\forall\; m\in \enne).
\]

\indent Thus, taking in $X$ the sequence $\{s_m\}$ which, for each $m\in \enne$, is given by $\displaystyle{(s_m)_k:=b_m}$ for every $k\in {\mathbb{Z}}[1,n]$,
owing to (\ref{immersione2}) and noting that
\begin{eqnarray}\label{11}
\Phi(u)&\leq& \frac{u^{t}Au}{2}+\sum_{k=1}^n\left(\int_0^{u_k}|h_k(t)|dt\right)\nonumber\\
 &\leq&\frac{u^{t}Au}{2}+\frac{L}{2}\sum_{k=1}^nu_k^2\nonumber\\
&=&\frac{u^{t}Au}{2}+\frac{L}{2}\|u\|^2.\nonumber
\end{eqnarray}
\noindent one immediately has
\begin{eqnarray*}\label{11}
 J_\lambda(s_m)&=& \frac{s_m^{t}As_m}{2}-{\lambda}\sum_{k=1}^{n}F_k(b_m)\nonumber\\
 &\leq& \frac{\trace +2\sum_{i<j}a_{ij}+nL}{2}b^2_m - \lambda\sum_{k=1}^{n}F_k(b_m)\nonumber\\
&<&\left(\frac{\trace +2\sum_{i<j}a_{ij}+nL}{2}-{\lambda}M \right)b^2_m.\nonumber
\end{eqnarray*}
that is, $\displaystyle\lim_{m\to \infty}J_\lambda(s_m)=-\infty$.

Next, assume that $B^\infty<+\infty$. Since $$\lambda >\displaystyle\frac{\trace +2\sum_{i<j}a_{ij}+nL}{2B^\infty},$$ we can fix $\varepsilon>0$ such that
$$ \displaystyle{\varepsilon<B^\infty-\frac{\trace +2\sum_{i<j}a_{ij}+nL}{2\lambda}}.$$\par
\indent  Therefore,
also calling $\{b_m\}$ a sequence of positive numbers such that $\displaystyle\lim_{m\rightarrow \infty}b_m=+\infty$ and
$$
(B^\infty-\varepsilon)b^2_m
<\sum_{k=1}^{n}F_k(b_m)<(B^\infty+\varepsilon) b^2_m, \quad \quad
(\forall\; m\in\enne)
$$
arguing as before and by choosing $\{s_m\}$ in $X$ as above, one has
$$J_\lambda(s_m)<\left(\frac{\trace +2\sum_{i<j}a_{ij}+nL}{2}-{\lambda}(B^\infty-\varepsilon) \right)b^2_m.$$
So, $\displaystyle\lim_{m \to \infty}J_\lambda(s_m)=-\infty$.

Hence, in both cases $J_\lambda$ is unbounded from below. The proof is thus complete.
\end{pf}

\begin{remark}\label{cp}
\rm{If $f_k$ are nonnegative continuous functions, condition $(\textrm{h}_\infty^L)$ reads as follows
\[
{\liminf_{t\rightarrow +\infty}}\frac{\displaystyle \sum_{k=1}^{n}F_k(t)}{t^2}< \displaystyle \frac{\lambda_1-L}{\trace +2\sum_{i<j}a_{ij}+nL}\limsup_{t\rightarrow+\infty}\frac{\displaystyle \sum_{k=1}^{n}F_k(t)}{t^2}.
\] }
\end{remark}

\indent Arguing as in the proof of Theorem \ref{Ma} and applying part $(\rm{b})$ of Theorem \ref{two}, we obtain the following result.
\begin{theorem}\label{Main2}
Let $A\in {\mathfrak{X}}_n$ and assume that the following inequality holds
\begin{itemize}
    \item [$(\textrm{h}_0^L)$] \qquad \qquad $A_0 <\displaystyle \frac{\lambda_1-L}{\trace +2\sum_{i<j}a_{ij}+nL}B^{0}.$
\end{itemize}

\noindent Then, for each
$$\lambda \in \displaystyle\left]
\frac{\trace +2\sum_{i<j}a_{ij}+nL}{2B^0},
\frac{\lambda_1-L}{2{{A}}_0}\right[,$$
 problem
$(S_{A,\lambda}^{f})$ admits a sequence of nontrivial solutions
$\{u_m\}$ such that $\displaystyle\lim_{m\rightarrow
\infty}\|u_m\|=\displaystyle\lim_{m\rightarrow
\infty}\|u_m\|_\infty=0$.
\end{theorem}
\section{Application}
\indent In this section we consider a discrete system, namely $(E_\lambda^{f,h})$, given as follows
$$
[u(i+1,j)-2u(i,j)+u(i-1,j)]+[u(i,j+1)-2u(i,j)+u(i,j-1)]
$$
$$
+\lambda f((i,j),u(i,j))+h(u(i,j))=0,\,\,\,\,\,\,\,\,\forall\,(i,j)\in {\mathbb{Z}}[1,m]\times {\mathbb{Z}}[1,n],
$$
with boundary conditions
$$
u(i,0)=u(i,n+1)=0,\,\,\,\,\,\,\forall\, i\in {\mathbb{Z}}[1,m],
$$
$$
u(0,j)=u(m+1,j)=0,\,\,\,\,\,\,\forall\, j\in {\mathbb{Z}}[1,n],
$$
where $f:{\mathbb{Z}}[1,m]\times{\mathbb{Z}}[1,n]\times \erre\rightarrow \erre$ denotes a continuous function, $\lambda$ is a positive real parameter and $h:\erre\rightarrow \erre$ be a Lipschitz continuous function with constant $L_h$.\par

As ambient space $X$, we consider the $mn$-dimensional Banach space $\erre^{mn}$
endowed by the norm
$$
\|u\|:=\Big(\sum_{k=1}^{mn}u_k^2\Big)^{1/2}.
$$

 Further, if $\ell\in\enne$, the symbol ${\mathfrak{M}}_{\ell\times \ell}(\erre)$ stands for the linear space of all the matrices of order $\ell$ with real entries.\par
  Let $v:{\mathbb{Z}}[1,m]\times{\mathbb{Z}}[1,n]\rightarrow {\mathbb{Z}}[1,mn]$ be the bijection defined by $v(i,j):=i+m(j-1),$ for every $(i,j)\in {\mathbb{Z}}[1,m]\times{\mathbb{Z}}[1,n]$.\par

 \indent Let us denote $w_k:=u({v^{-1}(k)})$ and $g_k(w_k):=f(v^{-1}(k),w_k)$, for every $k\in {\mathbb{Z}}[1,mn]$. With the above notations, problem $(E_\lambda^{f,h})$ can be written as a nonlinear algebraic system of the form
$$
Aw=\lambda g(w)+\widetilde{h}(w),\eqno\displaystyle{(S_{A,\lambda}^{g,\widetilde{h}})}
$$
where $A$ is given by
$$
A:=
\left(
  \begin{array}{ccccccccc}
    D & -I_m & 0 & 0 & ... & 0 & 0 & 0 & 0 \\
    -I_m & D & -I_m & 0 & ... & 0 & 0 &0 & 0 \\
    0 & -I_m & D & -I_m & ... & 0 & 0 & 0 & 0 \\
    0 & 0 & -I_m & D & ... & 0 & 0 & 0 & 0 \\
     &  &  &  & \ddots &  &  &  &  \\
    0 & 0 & 0 & 0 & ... & D & -I_m & 0 & 0 \\
    0 & 0 & 0 & 0 & ... & -I_m & D & -I_m & 0 \\
    0 & 0 & 0 & 0 & ... & 0 & -I_m & D & -I_m \\
    0 & 0 & 0 & 0 & ... & 0 & 0 & -I_m & D \\
  \end{array}
\right)\in {\mathfrak{M}}_{mn\times mn}(\erre),
$$
in which $D$ is defined by
$$
D:=
\left(
  \begin{array}{ccccccccc}
    4 & -1 & 0 & 0 & ... & 0 & 0 & 0 & 0 \\
    -1 & 4 & -1 & 0 & ... & 0 & 0 &0 & 0 \\
    0 & -1 & 4 & -1 & ... & 0 & 0 & 0 & 0 \\
    0 & 0 & -1 & 4 & ... & 0 & 0 & 0 & 0 \\
     &  &  &  & \ddots &  &  &  &  \\
    0 & 0 & 0 & 0 & ... & 4 & -1 & 0 & 0 \\
    0 & 0 & 0 & 0 & ... & -1 & 4 & -1 & 0 \\
    0 & 0 & 0 & 0 & ... & 0 & -1 & 4 & -1 \\
    0 & 0 & 0 & 0 & ... & 0 & 0 & -1 & 4 \\
  \end{array}
\right)\in {\mathfrak{M}}_{m\times m}(\erre),
$$
\noindent $I_m\in {\mathfrak{M}}_{m\times m}(\erre)$ is the identity matrix and $g(w):=(g_1(w_1),...,g_{mn}(w_{mn}))^t$, $\widetilde{h}(w):=(h(w_1),...,h(w_{mn}))^t$, for every $w\in X$.\par
 \indent In \cite{JY}, Ji and Yang studied the structure of the spectrum of the above (non-perturbed) Dirichlet problem. By their result we have that $A\in {\mathfrak{X}}_{mn}$.\par
 \indent It is easy to observe that the solutions of $(E_\lambda^{f,h})$ are the critical points of the $C^1$-functional
$$
J_\lambda(w):=\frac{w^{t}Aw}{2}-\lambda\sum_{k=1}^{mn}\int_0^{w_k}g_k(t)dt-\sum_{k=1}^{mn}\int_0^{w_k}h(t)dt,\,\,\,\,\,\forall\; w\in X.
$$

Denote by $\lambda_A$ the first eigenvalue of the matrix $A$. By using the above variational framework, Theorem \ref{Ma} assumes the following form.

\begin{theorem}\label{Ma2}
Assume that
$
\lambda_A<L_h,
$
in addition to
\begin{itemize}
    \item [$(\textrm{h}_\infty^{h})$] \qquad $\displaystyle\liminf_{t\rightarrow +\infty}\frac{\displaystyle{\sum_{k=1}^{mn}}\displaystyle{\max_{|\xi|\leq
t}\displaystyle\int_{0}^{\xi}g_k(s)ds}}{t^2} <\displaystyle \frac{\lambda_A-L_h}{(2+L_h)(m+n)}\limsup_{t
\rightarrow+\infty}\frac{\displaystyle \sum_{k=1}^{mn}\displaystyle\int_{0}^{t}g_k(s)ds}{t^2}.$
\end{itemize}

\noindent Then for each $$\lambda \in \displaystyle\left]
\frac{(2+L_h)(m+n)}{2B^\infty},
\frac{\lambda_A-L_h}{2{{A}}_\infty}\right[,$$ problem
$(E_\lambda^{f,h})$ admits an unbounded sequence of solutions.
\end{theorem}

\begin{remark}{\em Substituting
$\xi\to+\infty$ with $\xi\to0^+$ in Theorem \ref{Ma2},
the same statement as
Theorem \ref{Main2} is easily proved.
}
\end{remark}

\begin{remark}\label{din}
{\em We just point out that Theorem \ref{intro} in Introduction directly follows by Theorem \ref{Ma2} assuming that $L_h<\lambda_A$.
}
\end{remark}

\begin{remark}\rm{
We refer to the paper of Galewski and Orpel \cite{GO} for some multiplicity results on discrete partial difference equations as well as to the monograph of Cheng \cite{Cheng} for their discrete geometrical interpretation. See also the papers \cite{MarcuMolica, MolicaRepovs1, MolicaRepovs2,MolicaRepovs3} for recent contribution on discrete problems.}
\end{remark}
{\bf Acknowledgements.}  This paper was written when the first author was visiting professor at the University of Ljubljana in 2013. He expresses his gratitude to the host institution for warm hospitality.
  The manuscript was realized within the auspices of the GNAMPA Project 2013 titled {\it Problemi non-locali di tipo Laplaciano frazionario} and the SRA grants P1-0292-0101 and J1-5435-0101.

\end{document}